\documentclass[smallextended]{svjour3}
\usepackage{eurosym}
\usepackage{amsfonts}
\usepackage{graphicx}
\usepackage{fullpage}
\usepackage{subfigure}
\usepackage{amssymb}
\usepackage{amsmath}

\RequirePackage{fix-cm}
\smartqed

\begin{document}

\title{A new reproducing kernel approach for nonlinear fractional
three-point boundary value problems}
\author{Mehmet Giyas Sakar \and Onur Sald{\i}r }
\institute{Mehmet Giyas Sakar \at
	Yuzuncu Yil University, Faculty of Sciences, Department of
	Mathematics, 65080, Van, Turkey \\
	Tel.: +90 (432) 2251701-27895\\
	\email{giyassakar@hotmail.com}           
	\and
	Onur Sald{\i}r \at
	Yuzuncu Yil University, Faculty of Sciences, Department of
	Mathematics, 65080, Van, Turkey\\
	\email{onursaldir@gmail.com}           
} 
\date{Received: date / Accepted: date}
\maketitle

\begin{abstract}
In this article, a new reproducing kernel approach is developed for
obtaining numerical solution of nonlinear three-point boundary value
problems with fractional order. This approach is based on reproducing kernel
which is constructed by shifted Legendre polynomials. In considered problem,
fractional derivatives with respect to $\alpha$ and $\beta$ are defined in
Caputo sense. This method has been applied to some examples which have exact
solutions. In order to shows the robustness of the proposed
method, some numerical results are given in tabulated forms.

\end{abstract}







\section{Introduction}

In this paper, a new iterative reproducing kernel approach will be
constructed for obtaining the numerical solution of nonlinear fractional three-point
boundary value problem, 
\begin{eqnarray}  \label{eq1}
a_{2}(\xi){\ }^{c}D^{\alpha}z(\xi)+a_{1}(\xi){\ }^{c}D^{\beta}z(\xi)+a_{0}(%
\xi)z(\xi)=g(\xi,z(\xi),z^{\prime }(\xi)),\quad \xi\in\lbrack 0,1]
\end{eqnarray}
with following boundary conditions, 
\begin{eqnarray}  \label{eq2}
z(0)=\gamma_{0},\,z(\theta)=\gamma_{1},\ z(1)=\gamma_{2}, \,\,\ 0<\theta<1,\
1<\alpha\leq2,\,\ 0<\beta\leq1.
\end{eqnarray}
Here, $a_{0}(\xi),$ $a_{1}(\xi),$ $a_{2}(\xi)$ $\in$ $C^{2}(0,1)$ and $%
g(\xi,z)\in$ $L_{\rho}^{2}[0,1]$ are sufficiently smooth functions and fractional derivatives are taken in Caputo sense. Without
loss of generality, we pay regard to $z(0)=0$, $z(\theta)=0$ and $z(1)=0$.
Because, $z(0)=\gamma_{0}$, $z(\theta)=\gamma_{1}$ and $z(1)=\gamma_{2}$
boundary conditions can be easily reduced to $z(0)=0$, $z(\theta)=0$ and $%
z(1)=0$.

Nonlinear fractional multi-point boundary value problems appear in a
different area of applied mathematics and physics \cite{1,2,3,4,5,6,7} and
references therein. Many important studies have been concerned in
engineering and applied science such as dynamical systems, fluid mechanics,
control theory, oil industries, heat conduction can be well-turned by
fractional differential equations \cite{8,9,10}. Some applications,
qualitative behaviors of solution and numerical methods to find approximate
solution have been investigated for differential equation with fractional
order \cite{11,12,13,14}.

More particularly, it is not easy to directly get exact solutions to most
differential equations with fractional order. Hence, numerical techniques
are utilised largely. Actually, in recent times many efficient and
convenient methods have been developed such as the finite difference method 
\cite{15}, finite element method \cite{16}, homotopy perturbation method 
\cite{17}, Haar wavelet methods \cite{18}, Adomian decomposition method \cite%
{19}, collocation methods \cite{20}, homotopy analysis method \cite{21},
differential transform method \cite{22}, variational iteration method \cite%
{23}, reproducing kernel space method \cite{24,25} and so on \cite{26,27,28}.

In 1908, Zaremba firstly introduced reproducing kernel concept \cite{29}.
His resarches with regard to boundary value problems which includes
Dirichlet condition. Reproducing kernel method (RKM) produces a
solution in convergent series form for many differential, partial and
integro-differential equations. For more information, we refer to \cite%
{30,31}. Recently, this RKM is applied for different type of problem. For
example, fractional order nonlocal boundary value problems \cite{32},
Riccati differential equations \cite{33}, forced Duffing equations with
nonlocal boundary conditions \cite{34}, Bratu equations with fractional
order Caputo derivative \cite{35}, time-fractional Kawahara equation \cite%
{36}, two-point boundary value problem \cite{37}, nonlinear fractional
Volterra integro-differential equations \cite{38}.

Recently, Legendre reproducing kernel method is proposed for fractional
two-point boundary value problem of Bratu Type Equations \cite{39}. The main
motivation of this paper is to extend the Legendre reproducing kernel
approach for solving nonlinear three-point boundary value problem with
Caputo derivative.

The remainder part of the paper is prepared as follows: some fundamental
definitions of fractional calculus and the theory of reproducing kernel with
Legendre basis functions are given in Section 2. The structure of solution
with Legendre reproducing kernel is demonstrated in Section 3. In order to
show the effectiveness of the proposed method, some numerical findings are
reported in Section 4. Finally, the last section contains some conclusions.

\section{Preliminaries}

In this section, several significant concepts, definitions, theorems, and
properties are provided which will be used in this research. \newline
\newline
\textbf{Definition 2.1} Let $z(\xi)\in C[0,1]$ and $\xi\in[0,1]$. Then, the $%
\alpha$ order left Riemann-Liouville fractional integral operator is given
as \cite{8,12,13}: 
\begin{eqnarray*}
J_{0+}^{\alpha}z(\xi)=\frac{1}{\Gamma(\alpha)}\int\limits_{0}^{\xi} {%
(\xi-s)^{\alpha -1}z(s)ds},
\end{eqnarray*}
here $\Gamma(.)$ is Gamma function, $\alpha\geq0$ and $\xi>0$. \newline
\newline
\textbf{Definition 2.2} Let $z(\xi)\in AC[0,1]$ and $\xi\in[0,1]$. Then, the 
$\alpha$ order left Caputo differential operator is given as \cite{8,12,13}: 
\begin{eqnarray*}
{\ }^{c}D_{0+}^{\alpha}z(\xi)=\frac{1}{\Gamma(m-\alpha)}\int_{0}^{\xi} \frac{
\partial^{m}}{\partial \xi^{m}}\frac{z(s)}{(\xi-s)^{m-\alpha-1}}ds, \,\
m-1<\alpha< m, m\in\mathbb{N} \,\ \hbox{and} \,\ \xi>0.
\end{eqnarray*}%
\newline
\noindent\textbf{Definition 2.3} In order to construct polynomial type
reproducing kernel, the first kind shifted Legendre polynomials are defined
over the interval $[0,1]$. For obtaining these polynomials the following
iterative formula can be given: 
\begin{eqnarray*}
P_{0}(\xi) &=& 1, \\
P_{1}(\xi) &=& 2\xi-1, \\
&\vdots& \\
(n+1)P_{n+1}(\xi) &=&(2n+1)(2\xi-1)P_{n}(\xi)-nP_{n-1}(\xi), \,\ n=1,2,...
\end{eqnarray*}
The orthogonality requirement is 
\begin{eqnarray}
\langle P_{n},P_{m} \rangle=\int_{0}^{1}\rho_{[0,1]}(\xi)P_{n}(\xi)P_{m}
(\xi)d\xi=\left\{ 
\begin{array}{ll}
0, & n\neq m, \\ 
1, & n=m=0, \\ 
\frac{1}{2n+1}, & n=m\neq0,%
\end{array}
\right.
\end{eqnarray}
here, weighted function is taken as, 
\begin{eqnarray}  \label{eq4}
\rho_{[0,1]}(\xi)=1.
\end{eqnarray}
Legendre basis functions can be established so that this basis function
system satisfy the homogeneous boundary conditions as: 
\begin{eqnarray}  \label{eq5}
z(0)=0 \,\,\hbox{and}\,\ z(1)=0.
\end{eqnarray}
Eq. (\ref{eq5}) has a advantageous feature for solving boundary value
problems. Therefore, these basis functions for $j\geq2$ can be defined as; 
\begin{eqnarray}  \label{eq6}
\phi_{j}(\xi)= \left\{ 
\begin{array}{ll}
P_{j}(\xi)-P_{0}(\xi), & \hbox{$j$ is even,} \\ 
P_{j}(\xi)-P_{1}(\xi), & \hbox{$j$ is odd.}%
\end{array}
\right.
\end{eqnarray}
such that this system satisfy the conditions 
\begin{eqnarray}  \label{eq7}
\phi_{j}(0)=\phi_{j}(1)=0.
\end{eqnarray}
It is worth noting that the basis functions given in Eq. (\ref{eq6}) are
complete system. For more information about orthogonal polynomials, please
see \cite{41,42,43}.\newline

\noindent\textbf{Definition 2.4} Let $\Omega \neq \emptyset$, and $\mathbb{H}
$ with its inner product $\langle\cdot,\cdot\rangle_\mathbb{H}$ be a Hilbert
space of real-valued functions on $\Omega$. Then, the reproducing kernel of $%
\mathbb{H}$ is $R:\Omega\times \Omega\rightarrow \mathbb{R}$ iff

\begin{enumerate}
\item $R(\cdot,\xi) \in \mathbb{H}, \forall \xi \in \Omega$

\item $\langle\phi,R(\cdot,\xi) \rangle_\mathbb{H} = \phi(\xi),
\forall\phi\in \mathbb{H}, \forall \xi \in \Omega$.
\end{enumerate}

The last condition is known as reproducing property. Especially, for any $x$%
, $\xi$ $\in$ $\Omega$, 
\begin{eqnarray}
R(x,\xi)=\langle R(\cdot,x),R(\cdot,\xi) \rangle_\mathbb{H}.  \notag
\end{eqnarray}

If a Hilbert space satisfies the above two conditions then is called
reproducing kernel Hilbert space. Uniqueness of the reproducing kernel can
be shown by use of Riesz representation theorem \cite{40}. \newline
\newline
\textbf{Theorem 2.1} Let $\{e_{j}\}_{j=1}^{n}$ be an orthonormal basis of $n$%
-dimensional Hilbert space $\mathbb{H}$, then 
\begin{eqnarray}  \label{eq8}
R(x,\xi )=R_{x}(\xi )=\sum_{j=1}^{n}\bar{e}_{j}(x)e_{j}(\xi )
\end{eqnarray}
is reproducing kernel of $\mathbb{H}$ \cite{30,31}.\newline
\newline
\textbf{Definition 2.5} Let $W_{\rho }^{m}[0,1]$ polynomials space be
pre-Hilbert space over $[0,1]$ with real coefficients and its degree $\leq m$
and inner product as: 
\begin{equation}  \label{eq9}
\langle z,v\rangle _{W_{\rho }^{m}}=\int_{0}^{1}\rho _{\lbrack 0,1]}(\xi
)z(\xi )v(\xi )d\xi ,\,\,\ \forall z,v\in W_{\rho }^{m}[0,1]
\end{equation}
with $\rho _{\lbrack 0,1]}(\xi )$ described by Eq. (\ref{eq4}), and the norm 
\begin{equation}  \label{eq10}
\Vert z\Vert _{W_{\rho }^{m}}=\sqrt{\langle z,z\rangle }_{W_{\rho
}^{m}},\,\,\ \forall z\in W_{\rho }^{m}[0,1].
\end{equation}%
With the aid of definiton of $L^{2}$ Hilbert space, $L_{\rho
}^{2}[0,1]=\{g|\int_{0}^{1}\rho _{\lbrack 0,1]}(\xi )|g(\xi )|^{2}d\xi
<\infty \}$ for any fixed $m$, $W_{\rho }^{m}[0,1]$ is a subspace of $%
L_{\rho }^{2}[0,1]$ and $\forall z,v\in W_{\rho }^{m}[0,1]$, $\langle
z,v\rangle _{W_{\rho }^{m}}=\langle z,v\rangle _{L_{\rho }^{2}}$ \newline
\newline
\textbf{Theorem 2.2} $W_{\rho }^{m}[0,1]$ Hilbert space is a reproducing
kernel space. \newline
\newline
\textbf{Proof.} From Definition 2.5, it is quite apparent that $W_{\rho
}^{m}[0,1]$ functions space is a finite-dimensional. It is well known that
all finite-dimensional pre-Hilbert space is a Hilbert space. Herewith, using
this consequence and Theorem 2.1, $W_{\rho }^{m}[0,1]$ is a reproducing
kernel space.\newline
\newline
For solving problem (\ref{eq1})-(\ref{eq2}), it is required to describe a
closed subspace of $W_{\rho }^{m}[0,1]$ so that satisfy homogeneous boundary
conditions. \newline
\newline
\textbf{Definition 2.6} Let 
\begin{equation*}
{\ }^{0}W_{\rho }^{m}[0,1]=\{z\text{ }|\text{ }z\in W_{\rho }^{m}[0,1],\text{
}z(0)=z(1)=0\}.
\end{equation*}%
One can easily demonstrate that ${\ }^{0}W_{\rho }^{m}[0,1]$ is a
reproducing kernel space using Eq. (\ref{eq6}). From Theorem 2.1, the kernel
function $R_{x}^{m}(\xi )$ of ${\ \ }^{0}W_{\rho }^{m}[0,1]$ can be written
as 
\begin{equation}  \label{eq11}
R_{x}^{m}(\xi )=\sum_{j=2}^{m}h_{j}(\xi ){h_{j}}(x).
\end{equation}%
Here, $h_{j}(\xi )$ is complete system which is easily obtained from basis
functions in Eq. (\ref{eq6}) with the help of Gram-Schmidt
orthonormalization process. Eq. (\ref{eq11}) is very useful for
implementation. In other words, $R_{x}^{m}(\xi )$ and $W_{\rho }^{m}[0,1]$
can readily re-calculated by increasing $m$.

\section{Main Results}

In this section, some important results related to reproducing kernel method
with shifted Legendre polynomials are presented. In the first subsection,
generation of reproducing kernel which is satify three-point boundary value
problems is presented. In the second subsection, representation of solution
is given ${\ }^{\theta}W_{\rho}^{m}[0,1]$. Then, we will construct an
iterative process for nonlinear problem in third subsection.

\subsection{Generation of reproducing kernel for three-point boundary value
problems}

In this subsection, we shall generate a reproducing kernel Hilbert space ${\ 
}^{\theta }W_{\rho }^{m}[0,1]$ in which every functions satisfies $z(0)=0$, $%
z(\theta )=0$ and $z(1)=0$. \newline
\newline
${\ }^{\theta }W_{\rho }^{m}[0,1]$ is defined as ${\ }^{\theta }W_{\rho
}^{m}[0,1]=\{z|z\in W_{\rho }^{m}[0,1],z(0)=z(\theta )=z(1)=0\}$. \newline
\newline
Obviously, ${\ }^{\theta }W_{\rho }^{m}[0,1]$ reproducing kernel space is a
closed subspace of ${\ }^{0}W_{\rho }^{m}[0,1]$. The reproducing kernel of ${%
\ }^{\theta }W_{\rho }^{m}[0,1]$ can be given with the following theorem. 
\newline
\newline
\textbf{Theorem 3.1} The reproducing kernel ${\ }^{\theta }R_{x}^{m}(\xi )$
of ${\ }^{\theta }W_{\rho }^{m}[0,1]$, 
\begin{equation}  \label{eq12}
{\ }^{\theta }R_{x}^{m}(\xi )=R_{x}^{m}(\xi )-\frac{R_{x}^{m}(\theta
)R_{\theta }^{m}(\xi )}{R_{\theta }^{m}(\theta )}.
\end{equation}
\newline
\textbf{Proof.} Frankly, not all elements of ${\ }^{0}W_{\rho }^{m}[0,1]$
vanish at $\theta $. This shows that $R_{\theta }^{m}(\theta )\neq $ 0.
Hence, it can be easily seen that ${\ }^{\theta }R_{x}^{m}(\theta )={\ }%
^{\theta }R_{\theta }^{m}(\xi )=0$ and therefore ${\ }^{\theta
}R_{x}^{m}(\xi )\in {\ }^{\theta }W_{\rho }^{m}[0,1]$. For $\forall z\left(
x\right) \in $ $^{\theta }W_{\rho }^{m}[0,1]$, clearly, $z\left( \theta
\right) =0$, it follows that

\begin{equation*}
{<z(x),}^{\theta }R_{x}^{m}(\xi )>_{^{\theta }W_{\rho }^{m}[0,1]}={<z(x),}%
\text{ }R_{x}^{m}(\xi )>-\frac{R_{x}^{m}(\alpha )z(\theta )}{R_{\theta
}^{m}(\theta )}=z(\xi ).
\end{equation*}

Namely, ${\ }^{\theta }R_{x}^{m}(\xi )$ is of reproducing kernel of $%
^{\theta }W_{\rho }^{m}[0,1],$. This completes the proof.

\subsection{Representation of solution in ${\ }^{\protect\theta}W_{\protect%
\rho }^{m}[0,1]$ Hilbert space}

In this subsection, reproducing kernel method with Legendre polyomials is
established for obtaining numerical solution of three-point boundary value
problem. For Eqs. (\ref{eq1})-(\ref{eq2}), the approximate solution shall be
constructed in ${\ }^{\theta }W_{\rho }^{m}[0,1]$. Firstly, we will define
linear operator $L$ as follow, 
\begin{equation*}
L:{\ }^{\theta }W_{\rho }^{m}[0,1]\rightarrow L_{\rho }^{2}[0,1]
\end{equation*}%
such that 
\begin{equation*}
Lz(\xi ):=a_{2}(\xi ){\ }^{c}D^{\alpha }z(\xi )+a_{1}(\xi ){\ }^{c}D^{\beta
}z(\xi )+a_{0}(\xi )z(\xi ).
\end{equation*}%
The Eqs.(\ref{eq1})-(\ref{eq2}) can be stated as follows 
\begin{equation}
\left\{ 
\begin{array}{ll}
Lz=g(\xi ,z(\xi ),z^{\prime }(\xi )) &  \\ 
z(0)=z(\theta )=z(1)=0. & 
\end{array}%
\right.   \label{eq13}
\end{equation}%
Easily can be shown that linear operator $L$ is bounded. We will obtain the
representation solution of Eq. (\ref{eq13}) in the ${\ }^{\theta }W_{\rho
}^{m}[0,1]$ space. Let $^{\theta }R_{x}^{m}(\xi )$ be the polynomial form of
reproducing kernel in ${\ }^{\theta }W_{\rho }^{m}[0,1]$ space.\newline
\newline
\textbf{Theorem 3.2} Let $\{\xi _{j}\}_{j=0}^{m-2}$ be any $(m-1)$ distinct
points in open interval $(0,1)$ for Eqs. (\ref{eq1})-(\ref{eq2}), then $\psi
_{j}^{m}(\xi )=L^{\ast }$ $^{\theta }R_{\xi _{j}}^{m}(\xi )=L_{x}$ $^{\theta
}R_{x}^{m}(\xi )|_{x=\xi _{j}}.$\newline
\textbf{Proof.} For any fixed $\xi _{j}\in (0,1)$, put 
\begin{eqnarray}
\psi _{j}^{m}(\xi ) &=&L^{\ast \text{ }\theta }R_{\xi _{j}}^{m}(\xi
)=\langle L^{\ast \text{ }\theta }R_{\xi _{j}}^{m}(\xi ),^{\theta }R_{\xi
}^{m}(x)\rangle _{{\ }^{\theta }W_{\rho }^{m}}  \notag  \label{eq14} \\
&=&\langle ^{\theta }R_{\xi _{j}}^{m}(\xi ),L_{x}\text{ }^{\theta }R_{\xi
}^{m}(x)\rangle _{L_{\rho }^{2}}=L_{x}\text{ }^{\theta }R_{\xi
}^{m}(x)|_{x=\xi _{j}}.
\end{eqnarray}%
It is quite obvious that $^{\theta }R_{\xi }^{m}(x)=$ $^{\theta
}R_{x}^{m}(\xi )$. Therefore $\psi _{j}^{m}(\xi )=L^{\ast }$ $^{\theta
}R_{\xi _{j}}^{m}(\xi )=L_{x}$ $^{\theta }R_{x}^{m}(\xi )|_{x=\xi _{j}}$.
Here, $L^{\ast }$ shows the adjoint operator of $L$. For any fixed $m$ and $%
\xi _{j}\in (0,1)$, $\psi _{j}^{m}\in {\ }^{\theta }W_{\rho }^{m}[0,1]$.%
\newline
\newline
\textbf{Theorem 3.3} Let $\{\xi _{j}\}_{j=0}^{m-2}$ be any $(m-1)$ distinct
points in open interval $(0,1)$ for $m\geq 2$, then $\{\psi
_{j}^{m}\}_{j=0}^{m-2}$ is complete in ${\ }^{\theta }W_{\rho }^{m}[0,1]$.%
\newline
\newline
\textbf{Proof.} For every fixed $z\in {\ }^{\theta }W_{\rho }^{m}[0,1]$, let 
\begin{equation*}
\langle z(\xi ),\psi _{j}^{m}(\xi )\rangle _{{\ }^{\theta }W_{\rho }^{m}}=0,
\end{equation*}%
this result shows, for $j=0,1,...,m-2$, 
\begin{eqnarray}
\langle z(\xi ),\psi _{j}^{m}(\xi )\rangle _{{\ }^{\theta }W_{\rho }^{m}}
&=&\langle z(\xi ),L^{\ast \text{ }\theta }R_{\xi _{j}}^{m}(\xi )\rangle _{{%
\ }^{\theta }W_{\rho }^{m}}  \notag \\
&=&\langle Lz(\xi ),^{\theta }R_{\xi _{j}}^{m}(\xi )\rangle _{L_{\rho }^{2}}
\notag \\
&=&Lz(\xi _{j})=0.
\end{eqnarray}%
In Eq. (15), by use of inverse operator, it is decided that $z\equiv 0$.
Thus, $\{\psi _{j}^{m}\}_{j=0}^{m-2}$ is complete in ${\ }^{\theta }W_{\rho
}^{m}[0,1]$. This completes the proof. \newline
\newline
Theorem 3.3 indicates that in Legendre reproducing kernel approach, using a
finite distinct points are enough. But, in traditional reproducing kernel
method need to dense sqeuence on the interval. Namely, this new approach is
vary from traditional method in \cite{28,32,33,34,35,38}.\newline
\newline
The orthonormal system $\{\bar{\psi}_{j}^{m}\}_{j=0}^{m-2}$ of ${\ }^{\theta
}W_{\rho }^{m}[0,1]$ can be derived with the help of the Gram-Schmidt
orthogonalization process using $\{\psi _{j}^{m}\}_{j=0}^{m-2}$, 
\begin{equation}
\bar{\psi}_{j}^{m}(\xi )=\sum_{k=0}^{j}\beta _{jk}^{m}\psi _{k}^{m}(\xi ),
\label{eq16}
\end{equation}%
here $\beta _{jk}^{m}$ show the coefficients of orthogonalization. \newline
\newline
\textbf{Theorem 3.4} Suppose that $z_{m}$ is the exact solution of Eqs. (\ref%
{eq1})-(\ref{eq2}) and $\{\xi _{j}\}_{j=0}^{m-2}$ shows any $(m-1)$ distinct
points in open interval $(0,1)$, in that case 
\begin{equation}
z_{m}(\xi )=\sum_{j=0}^{m-2}\sum_{k=0}^{j}\beta _{jk}^{m}g(\xi
_{k},z_{m}(\xi _{k}), z_{m}^\prime(\xi _{k}))\bar{\psi}_{j}^{m}(\xi ).  \label{eq17}
\end{equation}%
\textbf{Proof.} Since $z_{m}\in {\ }^{\theta }W_{\rho }^{m}[0,1]$ from
Theorem 3.3 can be written 
\begin{equation*}
z_{m}(\xi )=\sum_{i=0}^{m-2}\langle z_{m}(\xi ),\bar{\psi}_{j}^{m}(\xi
)\rangle _{{\ }^{\theta }W_{\rho }^{m}}\bar{\psi}_{j}^{m}(\xi ).
\end{equation*}%
On the other part, using Eq. (\ref{eq14}) and Eq. (\ref{eq16}), we obtain $%
z_{m}(\xi )$ which is the precise solution of Eq. (\ref{eq10}) in ${\ }%
^{\theta }W_{\rho }^{m}[0,1]$ as,

\begin{eqnarray*}
z_{m}(\xi ) &=&\sum_{j=0}^{m-2}\langle z_{m}(\xi ),\bar{\psi}_{j}^{m}(\xi
)\rangle _{{\ }^{\theta }W_{\rho }^{m}}\bar{\psi}_{j}^{m}(\xi ) \\
&=&\sum_{j=0}^{m-2}\langle z_{m}(\xi ),\sum_{k=0}^{j}\beta _{jk}^{m}\psi
_{k}^{m}(\xi )\rangle _{{\ }^{\theta }W_{\rho }^{m}}\bar{\psi}_{j}^{m}(\xi )
\\
&=&\sum_{j=0}^{m-2}\sum_{k=0}^{j}\beta _{jk}^{m}\langle z_{m}(\xi ),\psi
_{k}^{m}(\xi )\rangle _{{\ }^{\theta }W_{\rho }^{m}}\bar{\psi}_{j}^{m}(\xi )
\\
&=&\sum_{j=0}^{m-2}\sum_{k=0}^{j}\beta _{jk}^{m}\langle z_{m}(\xi ),L^{\ast
}{}^{\theta }R_{\xi _{k}}^{m}(\xi )\rangle _{{\ }^{\theta }W_{\rho }^{m}}%
\bar{\psi}_{j}^{m}(\xi ) \\
&=&\sum_{j=0}^{m-2}\sum_{k=0}^{j}\beta _{jk}^{m}\langle Lz_{m}(\xi
),^{\theta }R_{\xi _{k}}^{m}(\xi )\rangle _{L_{\rho }^{2}}\bar{\psi}%
_{j}^{m}(\xi ) \\
&=&\sum_{j=0}^{m-2}\sum_{k=0}^{j}\beta _{jk}^{m}\langle g(\xi ,z_{m}(\xi
),z^\prime_{m}(\xi)),^{\theta }R_{\xi _{k}}^{m}(\xi )\rangle _{L_{\rho }^{2}}\bar{\psi}%
_{j}^{m}(\xi ) \\
&=&\sum_{j=0}^{m-2}\sum_{k=0}^{j}\beta _{jk}^{m}g(\xi _{k},z_{m}(\xi _{k}),z^\prime_{m}(\xi _{k}))\bar{\psi}_{j}^{m}(\xi ).
\end{eqnarray*}%
The proof is completed. \newline
\newline
\textbf{Theorem 3.5} If $z_{m}(\xi )\in {\ }^{\theta }W_{\rho }^{m}[0,1]$,
then $|z_{m}^{(s)}(\xi )|\leq F\Vert z_{m}\Vert _{{\ }^{\theta }W_{\rho
}^{m}}$ for $s=0,\ldots ,m-1$, where $F$ is a constant. \newline
\newline
\textbf{Proof.} We have $z_{m}^{(s)}\left( \xi \right) =\langle z_{m}\left(
x\right) ,\partial _{\xi }^{s}$ $^{\theta }R_{\xi }^{m}\left( x\right)
\rangle _{{\ }^{\theta }W_{\rho }^{m}}$ for any $\xi ,\,x\in \left[ {0,1}%
\right] $, $s=0,\ldots ,m-1.$ From the expression of $^{\theta }R_{\xi
}^{m}\left( x\right) $, it pursue that $\left\Vert {\partial _{\xi }^{s}}%
^{\theta }R_{\xi }^{m}\left( x\right) \right\Vert _{{\ }^{\theta }W_{\rho
}^{m}}\leq F_{s},\,s=0,\ldots ,m-1.$\newline
So, 
\begin{eqnarray*}
|z_{m}^{(s)}(\xi )| &=&|{\langle z_{m}(\xi ),\partial _{\xi }^{s}}\text{{\ }}%
^{\theta }R_{\xi }^{m}\left( x\right) {\rangle _{{\ }^{\theta }W_{\rho }^{m}}%
}| \\
&\leq &\Vert {z_{m}(\xi )}\Vert _{{\ }^{\theta }W_{\rho }^{m}[0,1]}\Vert {%
\partial _{\xi }^{s}}\text{{\ }}^{\theta }R_{\xi }^{m}\left( \xi \right)
\Vert _{{\ }^{\theta }W_{\rho }^{m}} \\
&\leq &F_{s}\Vert {z_{m}(\xi )}\Vert _{{\ }^{\theta }W_{\rho
}^{m}},s=0,\ldots ,m-1.
\end{eqnarray*}%
Therefore, $|z_{m}^{(s)}(\xi )|\leq \max \{F_{0},\ldots ,F_{m-1}\}\left\Vert 
{z_{m}\left( \xi \right) }\right\Vert _{{\ }^{\theta }W_{\rho
}^{m}},\,s=0,\ldots ,m-1$. \newline
\newline
\textbf{Theorem 3.6} $z_{m}(\xi )$ and its derivatives $z_{m}^{(s)}(\xi )$
are respectively uniformly converge to $z(\xi )$ and $z^{(s)}(\xi )$ ($%
s=0,\ldots ,m-1$). \newline
\newline
\textbf{Proof} By using Theorem 3.5 for any $\xi \in \lbrack 0,1]$ we get 
\begin{eqnarray*}
|z_{m}^{(s)}(\xi )-z^{(s)}(\xi )| &=&|\langle z_{m}(\xi )-z(\xi ),\partial
_{\xi }^{s}\text{ }^{\theta }R_{\xi }^{m}\left( \xi \right) \rangle |_{{\ }%
^{\theta }W_{\rho }^{m}} \\
&\leq &\Vert \partial _{\xi }^{s}\text{ }^{\theta }R_{\xi }^{m}\left( \xi
\right) \Vert _{{\ }^{\theta }W_{\rho }^{m}}\Vert z_{m}(\xi )-z(\xi )\Vert _{%
{\ }^{\theta }W_{\rho }^{m}} \\
&\leq &F_{s}\Vert z_{m}(\xi )-z(\xi )\Vert _{{\ }^{\theta }W_{\rho
}^{m}},\,\ s=0,\ldots ,m-1. 
\end{eqnarray*}%
where $F_{0},\ldots ,F_{m-1}$ are positive constants. Therefore, if $%
z_{m}(\xi )\rightarrow z(\xi )$ in the meaning of the norm of ${\ }^{\theta
}W_{\rho }^{m}[0,1]$ as $m\rightarrow \infty $, $z_{m}(\xi )$ and its
derivatives $z_{m}^{^{\prime }}(\xi ),\ldots ,z_{m}^{(m-1)}(\xi )$ are
respectively uniformly converge to $z(\xi )$ and its derivatives $%
z^{^{\prime }}(\xi ),\ldots ,z^{(m-1)}(\xi )$. This completes the proof.\\

If considered problem is linear, numerical solution can be directly get from (\ref{eq17}). But, for nonlinear problem the following iterative procedure can be construct. 
\subsection{Construction of iterative procedure}

In this subsection, we will use the following iterative sequence to overcome
the nonlinearity of the problem, $y_{m}(\xi )$, inserting, 
\begin{equation}  \label{eq18}
\left\{ {{\begin{array}{*{20}c} {Ly_{m,n}\left( \xi \right) = g\left(
{\xi,z_{m,n-1}(\xi),z^\prime_{m,n-1}(\xi)} \right)} \hfill \\ {z_{m,n}\left( \xi \right) = P_{m-1}
y_{m,n} (\xi)} \hfill \\ \end{array}}}\right.
\end{equation}
here, orthogonal projection operator is defined as $P_{m-1}:{\ }^{\theta
}W_{\rho }^{m}[0,1]\rightarrow span\{\bar{\psi}_{0}^{m},\bar{\psi}%
_{1}^{m},\ldots ,\bar{\psi}_{m-2}^{m}\}$ and $y_{m,n}(\xi )\in {\ }^{\theta
}W_{\rho }^{m}[0,1]$ shows the $n$-th iterative numerical solution of (\ref%
{eq18}). Then, the following important theorem will be given for iterative
procedure. \newline
\newline
\textbf{Theorem 3.7} If $\{\xi _{j}\}_{j=0}^{m-2}$ is distinct points in
open interval $(0,1)$, then 
\begin{equation}  \label{eq19}
y_{m,n}(\xi )=\sum_{j=0}^{m-2}\sum_{k=0}^{j}\beta _{jk}^{m}g(\xi
_{k},z_{m,n-1}(\xi _{k}),z^\prime_{m,n-1}(\xi _{k}))\bar{\psi}_{j}^{m}(\xi )
\end{equation}
\textbf{Proof.} Since $y_{m,n}(\xi )\in {\ }^{\theta }W_{\rho }^{m}[0,1]$, $%
\{\bar{\psi}_{j}^{m}(\xi )\}_{j=0}^{m-2}$ is the complete orthonormal system
in ${\ }^{\theta }W_{\rho }^{m}[0,1]$, 
\begin{eqnarray*}
y_{m,n}(\xi ) &=&\sum_{j=0}^{m-2}\langle y_{m,n}(\xi ),\bar{\psi}%
_{j}^{m}(\xi )\rangle _{{\ }^{\theta }W_{\rho }^{m}}\bar{\psi}_{j}^{m}(\xi )
\\
&=&\sum_{j=0}^{m-2}\langle y_{m,n}(\xi ),\sum_{k=0}^{j}\beta _{jk}^{m}\psi
_{k}^{m}(\xi )\rangle _{{\ }^{\theta }W_{\rho }^{m}}\bar{\psi}_{j}^{m}(\xi )
\\
&=&\sum_{j=0}^{m-2}\sum_{k=0}^{j}\beta _{jk}^{m}\langle y_{m,n}(\xi ),\psi
_{k}^{m}(\xi )\rangle _{{\ }^{\theta }W_{\rho }^{m}}\bar{\psi}_{j}^{m}(\xi )
\\
&=&\sum_{j=0}^{m-2}\sum_{k=0}^{j}\beta _{jk}^{m}\langle y_{m,n}(\xi
),L^{\ast \text{ }\theta }R_{\xi _{k}}^{m}(\xi )\rangle _{{\ }^{\theta
}W_{\rho }^{m}}\bar{\psi}_{j}^{m}(\xi ) \\
&=&\sum_{j=0}^{m-2}\sum_{k=0}^{j}\beta _{jk}^{m}\langle Ly_{m,n}(\xi
),^{\theta }R_{\xi _{k}}^{m}(\xi )\rangle _{L_{\rho }^{2}}\bar{\psi}%
_{j}^{m}(\xi ) \\
&=&\sum_{j=0}^{m-2}\sum_{k=0}^{j}\beta _{jk}^{m}\langle g(\xi ,z_{m,n-1}(\xi
),z^\prime_{m,n-1}(\xi
)),^{\theta }R_{\xi _{k}}^{m}(\xi )\rangle _{L_{\rho }^{2}}\bar{\psi}%
_{j}^{m}(\xi ) \\
&=&\sum_{j=0}^{m-2}\sum_{k=0}^{j}\beta _{jk}^{m}g(\xi _{k},z_{m,n-1}(\xi
_{k}),z^\prime_{m,n-1}(\xi
_{k}))\bar{\psi}_{j}^{m}(\xi )
\end{eqnarray*}%
This completes the proof.\newline
Taking $z_{m,0}(\xi )=0$ and define the iterative sequence 
\begin{equation}  \label{eq20}
z_{m,n}(\xi )=P_{m-1}y_{m,n}(\xi )=\sum_{j=0}^{m-2}\sum_{k=0}^{j}\beta
_{jk}^{m}g(\xi _{k},z_{m,n-1}(\xi _{k}),z^\prime_{m,n-1}(\xi _{k}))\bar{\psi}_{j}^{m}(\xi ),\,\
n=1,2,\ldots
\end{equation}

\section{Numerical applications}

In this section, some nonlinear three-point boundary value problems are
considered to exemplify the accuracy and efficiency of proposed approach.
Numerical results which is achieved by L-RKM are shown with tables.\newline
\newline
\textbf{Example 4.1} We consider the following fractional order nonlinear
three-point boundary value problem with Caputo derivative: 
\begin{equation}  \label{eq21}
{\ }^{c}D^{\alpha }z(\xi ) + (\xi +1) {\ }^{c}D^{\beta }z(\xi ) + \xi
z(\xi)-z^{2}(\xi )=f(\xi),\quad 1<\alpha \leq 2.\quad 0<\beta \leq 1.
\end{equation}
\begin{equation}  \label{eq22}
z(0)=z(\frac{1}{2}) =z(1)=0.
\end{equation}%
Here, $f(\xi)$ a known function such that the exact solution of this problem is $z(\xi )=\xi(\xi-\frac{1}{2})(\xi-1)$.

By using proposed approach for Eqs. (\ref{eq21})-(\ref{eq22}), and choosing
nodal points as $\xi _{j}=\frac{j+0.3}{m},\,j=0,1,\,2,...,m-2$, the
approximate solution $z_{m,n}\left( \xi \right)$ is computed by Eq. (\ref%
{eq20}). For (\ref{eq21})-(\ref{eq22}), comparison of absolute errors for
different $\alpha$, $\beta$ values are demonstrated in Table 1 and Table 2
and comparison of exact solution and numerical solution for $\alpha=1.75$
and $\beta=0.75$ is given in Table 3. \newline

\noindent \textbf{Example 4.2} We take care of the following nonlinear
three-point boundary value problem with Caputo derivative 
\begin{equation}  \label{eq23}
\xi^2{\ }^{c}D^{\alpha }z(\xi ) + (\xi^2-1) {\ }^{c}D^{\beta }z(\xi ) +
\xi^3 z(\xi)-z(\xi)z^\prime(\xi)-z^3(\xi)=f(\xi),\quad 1<\alpha \leq 2.\quad
0<\beta \leq 1.
\end{equation}
\begin{equation}  \label{eq24}
z(0)=z(\frac{3}{5}) =z(1)=0.
\end{equation}
Here, $f(\xi)$ a known function such that the exact solution of this problem is $z(\xi )=\xi(\xi-\frac{3}{5})(\xi-1)$.

By using proposed approach for Eqs. (\ref{eq23})-(\ref{eq24}), and choosing
nodal points as $\xi _{j}=\frac{j+0.3}{m},\,j=0,1,\,2,...,m-2$, the
approximate solution $z_{m,n}\left( \xi \right)$ is computed by Eq. (\ref%
{eq20}). For (\ref{eq23})-(\ref{eq24}), comparison of absolute errors for
different $\alpha$, $\beta$ values are demonstrated in Table 4 and Table 5
and comparison of exact solution and numerical solution for $\alpha=1.75$
and $\beta=0.75$ is given in Table 6. \newline

\section{Conclusion}

In this research, a novel numerical approach which is called L-RKM has been proposed and
successfully implemented to find the approximate solution of nonlinear
three-point boundary value problems with Caputo derivative. For nonlinear
problem, a new iterative process is proposed. Numerical findings show that
the present approach is efficient and convenient for solving three-point
boundary value problems with fractional order.

\section*{Tables}

\begin{table}[h!]
\caption{Comparison absolute error of Example 4.1 for various $\protect\alpha%
, \protect\beta$ ($m=3$, $n=3$)}%
\begin{tabular}{cccccc}
\hline
$x$ & $\alpha=2, \beta=1 $ & $\alpha=1.9, \beta=0.9$ & $\alpha=1.8, \beta
=0.8$ & $\alpha=1.7, \beta= 0.7$ & $\alpha=1.6, \beta= 0.6$ \\ \hline
0.1 & 3.25E-11 & 6.08E-11 & 3.70E-12 & 4.45E-10 & 3.14E-9 \\ 
0.2 & 5.45E-11 & 8.75E-11 & 4.93E-12 & 8.70E-10 & 5.51E-9 \\ 
0.3 & 6.97E-11 & 9.00E-11 & 4.34E-12 & 1.32E-9 & 7.36E-9 \\ 
0.4 & 8.17E-11 & 7.81E-11 & 2.55E-12 & 1.85E-9 & 8.96E-9 \\ 
0.5 & 0 & 0 & 0 & 0 & 0 \\ 
0.6 & 1.10E-10 & 5.14E-11 & 2.03E-12 & 3.36E-9 & 1.24E-8 \\ 
0.7 & 1.34E-10 & 5.65E-11 & 3.56E-12 & 4.44E-9 & 1.49E-8 \\ 
0.8 & 1.70E-10 & 8.70E-11 & 3.73E-12 & 5.80E-9 & 1.81E-8 \\ 
0.9 & 2.21E-10 & 1.53E-11 & 1.89E-12 & 7.49E-9 & 2.25E-8 \\ \hline
\end{tabular}%
\end{table}

\begin{table}[h!]
\caption{Comparison absolute error of Example 4.1 for various $\protect\alpha%
, \protect\beta$ ($m=3$, $n=5$)}%
\begin{tabular}{cccccc}
\hline
$x$ & $\alpha=2, \beta=1 $ & $\alpha=1.9, \beta=0.9$ & $\alpha=1.8, \beta
=0.8$ & $\alpha=1.7, \beta= 0.7$ & $\alpha=1.6, \beta= 0.6$ \\ \hline
0.1 & 3.78E-17 & 1.33E-16 & 2.49E-19 & 4.86E-15 & 7.48E-14 \\ 
0.2 & 5.27E-17 & 1.94E-16 & 3.36E-19 & 1.13E-14 & 1.36E-13 \\ 
0.3 & 5.10E-17 & 2.03E-16 & 3.80E-19 & 1.99E-14 & 1.91E-13 \\ 
0.4 & 3.91E-17 & 1.81E-16 & 5.20E-19 & 3.09E-14 & 2.43E-13 \\ 
0.5 & 0 & 0 & 0 & 0 & 0 \\ 
0.6 & 1.05E-17 & 1.36E-16 & 1.56E-18 & 6.18E-14 & 3.60E-13 \\ 
0.7 & 6.65E-18 & 1.57E-16 & 2.74E-18 & 8.25E-14 & 4.36E-13 \\ 
0.8 & 1.81E-17 & 2.34E-16 & 4.50E-18 & 1.07E-14 & 5.30E-13 \\ 
0.9 & 5.14E-17 & 3.91E-16 & 7.02E-18 & 1.36E-14 & 6.48E-13 \\ \hline
\end{tabular}%
\end{table}

\begin{table}[h!]
\caption{Numerical results of Example 4.1 for $m=5$, $n=9$ values ($\protect%
\alpha=1.75$, $\protect\beta=0.75$)}%
\begin{tabular}{cccc}
\hline
$x$ & Exact Sol. & Approximate Sol. & Absolute Error \\ \hline
0.0 & 0.000000000000000000000 & 0.000000000000000000000 & 0 \\ 
0.1 & 0.036000000000000000000 & 0.036000000000000000018 & 1.80E-20 \\ 
0.2 & 0.048000000000000000000 & 0.048000000000000000044 & 4.40E-20 \\ 
0.3 & 0.042000000000000000000 & 0.042000000000000000061 & 6.10E-20 \\ 
0.4 & 0.024000000000000000000 & 0.024000000000000000071 & 7.10E-20 \\ 
0.5 & 0.000000000000000000000 & 0.000000000000000000000 & 0 \\ 
0.6 & -0.024000000000000000000 & -0.023999999999999999899 & 1.01E-19 \\ 
0.7 & -0.042000000000000000000 & -0.041999999999999999819 & 1.81E-19 \\ 
0.8 & -0.048000000000000000000 & -0.047999999999999999694 & 3.06E-19 \\ 
0.9 & -0.036000000000000000000 & -0.035999999999999999491 & 5.09E-19 \\ 
1.0 & 0.000000000000000000000 & 0.000000000000000000000 & 0 \\ \hline
\end{tabular}%
\end{table}

\begin{table}[h!]
\caption{Comparison absolute error of Example 4.2 for various $\protect\alpha%
, \protect\beta$ ($m=3$, $n=8$)}
\begin{tabular}{cccccc}
\hline
$x$ & $\alpha=2, \beta=1 $ & $\alpha=1.9, \beta=0.9$ & $\alpha=1.8, \beta
=0.8$ & $\alpha=1.7, \beta= 0.7$ & $\alpha=1.6, \beta= 0.6$ \\ \hline
0.1 & 5.00E-11 & 4.11E-15 & 2.78E-13 & 6.10E-12 & 1.40E-12 \\ 
0.2 & 8.39E-11 & 3.45E-15 & 1.87E-13 & 1.17E-11 & 2.42E-11 \\ 
0.3 & 1.06E-10 & 1.24E-15 & 1.87E-13 & 1.45E-11 & 7.08E-11 \\ 
0.4 & 1.23E-10 & 9.21E-15 & 7.60E-13 & 1.20E-11 & 1.32E-10 \\ 
0.5 & 1.39E-10 & 1.97E-14 & 1.44E-12 & 1.75E-12 & 2.02E-10 \\ 
0.6 & 0 & 0 & 0 & 0 & 0 \\ 
0.7 & 1.88E-10 & 4.52E-14 & 2.81E-12 & 5.17E-11 & 3.45E-10 \\ 
0.8 & 2.31E-10 & 5.88E-14 & 3.32E-12 & 9.98E-11 & 4.05E-10 \\ 
0.9 & 2.93E-10 & 7.18E-14 & 3.60E-12 & 1.65E-10 & 4.49E-10 \\ \hline
\end{tabular}%
\end{table}

\begin{table}[h!]
\caption{Comparison absolute error of Example 4.2 for various $\protect\alpha%
, \protect\beta$ ($m=3$, $n=10$)}
\begin{tabular}{cccccc}
\hline
$x$ & $\alpha=2, \beta=1 $ & $\alpha=1.9, \beta=0.9$ & $\alpha=1.8, \beta
=0.8$ & $\alpha=1.7, \beta= 0.7$ & $\alpha=1.6, \beta= 0.6$ \\ \hline
0.1 & 3.49E-13 & 1.93E-18 & 4.48E-16 & 2.68E-14 & 4.75E-15 \\ 
0.2 & 5.87E-13 & 8.32E-19 & 2.99E-16 & 5.58E-14 & 7.68E-14 \\ 
0.3 & 7.47E-13 & 2.90E-18 & 3.06E-16 & 7.18E-14 & 2.25E-13 \\ 
0.4 & 8.65E-13 & 8.81E-18 & 1.23E-15 & 6.00E-14 & 4.21E-13 \\ 
0.5 & 9.76E-13 & 1.65E-17 & 2.34E-15 & 5.54E-15 & 6.45E-13 \\ 
0.6 & 0 & 0 & 0 & 0 & 0 \\ 
0.7 & 1.31E-12 & 3.55E-17 & 4.55E-15 & 2.91E-13 & 1.09E-12 \\ 
0.8 & 1.61E-12 & 4.60E-17 & 5.37E-15 & 5.64E-13 & 1.29E-12 \\ 
0.9 & 2.05E-12 & 5.65E-17 & 5.83E-15 & 9.39E-13 & 1.43E-12 \\ \hline
\end{tabular}%
\end{table}

\begin{table}[h!]
\caption{Numerical results of Example 4.2 for $m=5$, $n=9$ values ($\protect%
\alpha=1.75$ , $\protect\beta=0.75$)}
\begin{tabular}{cccc}
\hline
$x$ & Exact Sol. & Approximate Sol. & Absolute Error \\ \hline
0.0 & 0.000000000000000000000 & 0.000000000000000000000 & 0 \\ 
0.1 & 0.045000000000000000000 & 0.045000000000480782793 & 4.80E-13 \\ 
0.2 & 0.064000000000000000000 & 0.064000000000580045412 & 5.80E-13 \\ 
0.3 & 0.063000000000000000000 & 0.063000000000488602930 & 4.88E-13 \\ 
0.4 & 0.048000000000000000000 & 0.048000000000398645450 & 3.98E-13 \\ 
0.5 & 0.025000000000000000000 & 0.025000000000468872800 & 4.68E-13 \\ 
0.6 & 0.000000000000000000000 & 0.000000000000000000000 & 0 \\ 
0.7 & -0.021000000000000000000 & -0.020999999998651962040 & 1.34E-12 \\ 
0.8 & -0.032000000000000000000 & -0.031999999998006864020 & 1.99E-12 \\ 
0.9 & -0.027000000000000000000 & -0.026999999997598991320 & 2.40E-12 \\ 
1.0 & 0.000000000000000000000 & 0.000000000000000000000 & 0 \\ \hline
\end{tabular}%
\end{table}

\pagebreak


\begin{thebibliography}{99}
\bibitem{1} Y. Lin, J. Niu, M. Cui. A numerical solution to nonlinear second
order three-point boundary value problems in the reproducing kernel space.
Applied Mathematics and Computation, 218(14) (2012) 7362-7368.

\bibitem{2} M. Rehman, R. A. Khan, N. A. Asif. Three point boundary value
problems for nonlinear fractional differential equations. Acta Mathematica
Scientia, 31 (4) (2011) 1337-1346.

\bibitem{3} F. Geng, Solving singular second order three-point boundary
value problems using reproducing kernel Hilbert space method. Applied
Mathematics and Computation (2009) 215(6) 2095-2102.

\bibitem{4} C. P. Zhang, J. Niu, Y. Z. Lin. Numerical solutions for the
three-point boundary value problem of nonlinear fractional differential
equations. Abstract and Applied Analysis (2012) Volume 2012, Article ID
360631, 16 pages.

\bibitem{5} S. Etemad, S. K. Ntouyas, J. Tariboon. Existence results for
three-point boundary value problems for nonlinear fractional differential
equations. J. Nonlinear Sci. Appl. 9 (2016) 2105-2116.

\bibitem{6} B. Wu, X. Li. Application of reproducing kernel method to third
order three-point boundary value problems. Applied Mathematics and
Computation (2010) 217 (7) 3425-3428.

\bibitem{7} B. Ahmad, M. Alghanmi, S. K. Ntouyas, A. Alsaedi. A study of
fractional differential equations and inclusions involving generalized
Caputo-type derivative equipped with generalized fractional integral
boundary conditions. AIMS Mathematics, 4 (1) (2018) 26-42.

\bibitem{8} I. Podlubny. Fractional differential equations. Academic Press,
New York, 1999.

\bibitem{9} V. Lakshmikantham, S. Leela, J. Vasundhara Devi, Theory of
fractional dynamic systems. Cambridge Scientific Publishers, 2009.

\bibitem{10} R. Hilfer, Applications of Fractional Calculus in Physics,
World Scientific, Singapore, 2000

\bibitem{11} V. E. Tarasov, Fractional Dynamics: Application of Fractional
Calculus to Dynamics of Particles, Fields and Media, Springer, HEP, 2011.

\bibitem{12} K. Diethelm, The analysis of fractional differential equations.
Lecture notes in mathematics. Berlin Heidelberg: Springer-Verlag, 2010.

\bibitem{13} A. A. Kilbas, H. M. Srivastava, J. J. Trujillo, Theory and
Applications of Fractional Differential Equations, B.V: Elsevier Science,
2006.

\bibitem{14} A. Khalouta, A. Kadem. A new numerical technique for solving
Caputo time-fractional biological population equation. AIMS Mathematics, 4
(5) (2019) 1307-1319.

\bibitem{15} C. Tadjeran, M. M. Meerschaert. A second-order accurate
numerical method for the two-dimensional fractional diffusion equation. J.
Comput. Phys. 220 (2007) 813-823.

\bibitem{16} A. Esen, O. Tasbozan. Numerical solution of time fractional
Burgers equation by cubic B-spline finite elements. Mediterr. J. Math. 13
(2016) 1325-1337.

\bibitem{17} M. G. Sakar, F. Uludag, F. Erdogan, Numerical solution of
time-fractional nonlinear PDEs with proportional delays by homotopy
perturbation method. Appl. Math. Model. 40 (2016), 6639-6649.

\bibitem{18} U. Saeed, M. Rehman, Haar wavelet-quasilinearization technique
for fractional nonlinear differential equations. Applied Mathematics and
Computation 220 (2013) 630-648.

\bibitem{19} E. Babolian, A. R. Vahidi, A. Shoja, An efficient method for
nonlinear fractional differential equations: combination of the Adomian
decomposition method and spectral method, Indian J. Pure Appl. Math. (2014)
1017-1028.

\bibitem{20} L. Pezza, F. Pitolli, A multiscale collocation method for
fractional differential problems. Mathematics and Computers in Simulation
147 (2018) 210-219.

\bibitem{21} M. G. Sakar, F. Erdogan, The homotopy analysis method for
solving the time-fractional Fornberg-Whitham equation and comparison with
Adomian's decomposition method, Appl. Math. Model. 37 (20-21) (2013)
1634-1641.

\bibitem{22} H. Jafari, H. K. Jassim, S. P. Moshokoa, V. M. Ariyan, F.
Tchier, Reduced differential transform method for partial differential
equations within local fractional derivative operators, Advances in
Mechanical Engineering, 8 (4) (2016) 1\^{a}\euro ``6.

\bibitem{23} M. G. Sakar, O. Sald{\i}r, Improving variational iteration
method with auxiliary parameter for nonlinear time-fractional partial
differential equations. Journal of Optimization Theory and Applications 174
(2) (2017) 530-549.

\bibitem{24} M. Q. Xu, Y. Z. Lin, Simplified reproducing kernel method for
fractional differential equations with delay, Applied Mathematics Letters 52
(2016) 156-161.

\bibitem{25} Y. L. Wang, M. J. Du, C. L. Temuer, D. Tian, Using reproducing
kernel for solving a class of time-fractional telegraph equation with
initial value conditions, International Journal of Computer Mathematics 95
(8) (2018) 1609-1621.

\bibitem{26} A. Kadem, D. Baleanu. Fractional radiative transfer equation
within Chebyshev spectral approach. Computers and Mathematics with
Applications 59 (2010) 1865-1873.

\bibitem{27} S. S. E. Eldien, R. M. Hafez, A. H. Bhrawy, D. Baleanu, A. A.
E.Kalaawy. New numerical approach for fractional variational problems using
shifted Legendre orthonormal polynomials. Journal of Optimization Theory and
Applications 174 (2017) 295-320.

\bibitem{28} M. G. Sakar, A. Akg\"{u}l, D. Baleanu, On solutions of
fractional Riccati differential equations, Advances in Difference Equations
(2017) 2017:39.

\bibitem{29} S. Zaremba, Sur le calcul num\'{e}rique des fonctions demand%
\'{e}es dans le probl\'{e}me de Dirichlet et le probl\`{e}me hydrodynamique.
Bulletin International de l'Acad\'{e}mie des Sciences de Cracovie (1908) pp.
125-195.

\bibitem{30} M. Cui, Y. Lin. Nonlinear Numerical Analysis in the Reproducing
Kernel Space. Nova Science, New York, 2009.

\bibitem{31} S. Saitoh, Y. Sawano. Theory of Reproducing Kernels and
Applications, Springer, Singapore, 2016.

\bibitem{32} F. Geng, M. Cui, A reproducing kernel method for solving
nonlocal fractional boundary value problems. Applied Mathematics Letters  25
(5) (2012) 818-823.

\bibitem{33} M. G. Sakar, Iterative reproducing kernel Hilbert spaces method
for Riccati differential equation. Journal of Computational and Applied
Mathematics. 309 (2017) 163-174.

\bibitem{34} F. Geng, M. Cui, New method based on the HPM and RKHSM for
solving forced Duffing equations with integral boundary conditions. Journal
of Computational and Applied Mathematics. 233 (2009) 165-172.

\bibitem{35} E. Babolian, S. Javadi, E. Moradi, RKM for solving Bratu-type
differential equations of fractional order. Mathematical Methods in the
Applied Sciences. 39 (6) (2016) 1548-1557.

\bibitem{36} O. Sald{\i}r, M. G. Sakar, F. Erdogan. Numerical solution of
time-fractional Kawahara equation using reproducing kernel method with error
estimate. Computational and Applied Mathematics, 38 (4) (2019) 98.

\bibitem{37} M. Khaleghi, E. Babolian, S. Abbasbandy. Chebyshev reproducing
kernel method: application to two-point boundary value problems. Advances in
Difference Equations (2017) 26. DOI: 10.1186/s13662-017-1089-2

\bibitem{38} W. Jiang, T. Tian. Numerical solution of nonlinear Volterra
integro-differential equations of fractional order by the reproducing kernel
method. Applied Mathematical Modelling 39 (16) (2015) 4871-4876.

\bibitem{39} M. G. Sakar, O. Sald{\i}r, A. Akg\"{u}l. Numerical solution of
fractional Bratu type equations with Legendre reproducing kernel method.
International Journal of Applied and Computational Mathematics (2018) 4:126.

\bibitem{40} N. Aronszajn. Theory of reproducing kernels. Trans. Am. Math.
Soc. (1950), 68:337-404.

\bibitem{41} W. Kaplan, Advanced Calculus (5E). Pearson Education 2002.

\bibitem{42} E. D. Rainville. Special Functions. Chelsea Publishing Co., New
York (1960).

\bibitem{43} G. Szeg\"{o}. Orthogonal Polynomials. American Mathematical
Society Colloquium Publications, Providence, Rhode Island, 1939.

\newpage
\end{thebibliography}
\end{document}